\title{Unit Distances in Three Dimensions\thanks{%
  Work by Haim Kaplan has been supported
  by Grant 2006/204 from the U.S.-Israel Binational Science Foundation,
  and by grant 822/10 from the Israel Science Fund.
Work by Ji\v{r}\'{\i} Matou\v{s}ek  has been
partially supported by the  ERC Advanced Grant No.~267165.
Work by Zuzana Safernov\'a has been supported by the Charles University
grant GAUK~421511.
  Work by Micha Sharir has been supported
  by NSF Grant CCF-08-30272,
  by Grant 2006/194 from the U.S.-Israel Binational Science Foundation,
  by Grant 338/09 from the Israel Science Fund,
  and by the Hermann Minkowski--MINERVA Center for Geometry at Tel Aviv
  University.}
}

\documentclass[11pt]{article}

\usepackage{latexsym,color,graphicx,amsmath,amssymb}

\def\dg{{\rm deg}}
\def\Deg{{D}}
\def\Degg{{E}}

\def\reals{{\mathbb R}}
\def\eps{{\varepsilon}}
\newcommand\R{\reals}
\newcommand{\heading}[1]{\vspace{1ex}\par\noindent{\bf\boldmath #1}}

\newtheorem{theorem}{Theorem}
\newtheorem{lemma}[theorem]{Lemma}

\date{Rev. 5/VII/11 JM}

\author{
{\sc Haim Kaplan}
\\
    {\footnotesize School of Computer Science, }\\[-1.5mm]
    {\footnotesize Tel Aviv University, }\\[-1.5mm]
    {\footnotesize Tel~Aviv 69978, Israel}
\and {\sc Ji\v{r}\'{\i} Matou\v{s}ek}
\\
   {\footnotesize Department of Applied Mathematics and}\\[-1.5mm]
   {\footnotesize Institute of Theoretical Computer Science (ITI)}\\[-1.5mm]
   {\footnotesize  Charles University, Malostransk\'{e} n\'{a}m. 25}\\[-1.5mm]
   {\footnotesize  118~00~~Praha~1, Czech Republic, and}\\%[-1.5mm]
   {\footnotesize    Institute of  Theoretical Computer Science}\\[-1.5mm]
   {\footnotesize    ETH Zurich, 8092 Zurich, Switzerland}
\and
{\sc Zuzana Safernov\'a}
\\
{\footnotesize Department of Applied Mathematics}\\[-1.5mm]
   {\footnotesize  Charles University, Malostransk\'{e} n\'{a}m. 25}\\[-1.5mm]
   {\footnotesize  118~00~~Praha~1, Czech Republic}
\and {\sc Micha Sharir}
\\
    {\footnotesize School of Computer Science, }\\[-1.5mm]
    {\footnotesize Tel Aviv University, }\\[-1.5mm]
    {\footnotesize Tel~Aviv 69978, Israel, and}\\%[-1.5mm]
    {\footnotesize Courant Institute of Mathematical Sciences, }\\[-1.5mm]
    {\footnotesize New York University, }\\[-1.5mm]
    {\footnotesize New York, NY~~10012,~USA}
}

\begin{document}

\maketitle

\begin{abstract}
We show that the number of unit distances determined by $n$ points in
$\reals^3$ is $O(n^{3/2})$, slightly improving the bound of
Clarkson et al.~\cite{CEGSW}, established in 1990. The new proof uses
the recently introduced polynomial partitioning technique of Guth and
Katz~\cite{GK2}. While this paper was still in a draft stage,
a similar proof of our main result was posted to the arXiv
by Joshua Zahl~\cite{Za}.
\end{abstract}

\section{Introduction}

Let $P$ be a set of $n$ points in Euclidean $d$-dimensional space
$\reals^d$. What is the maximum possible number
of pairs of points in $P$ with distance exactly~$1$?
A standard construction, attributed to Lenz~\cite{Le55}, shows that
this number can be $\Theta(n^2)$ in $d\ge 4$ dimensions, so the only
interesting cases are $d=2,3$. The planar version is the classical
\emph{unit distances} problem of Erd{\H o}s~\cite{E46}, posed in 1946,
for which we refer to the literature (in particular, see
\cite{CEGSW,SST,Sz,Va}).
Here we focus on the case $d=3$.  This has been studied, back in 1990,
by Clarkson et al.~\cite{CEGSW}, who have established the
upper bound $O(n^{3/2}2^{O(\alpha^2(n))})$, where
$\alpha(\cdot)$ is the inverse Ackermann function.

In this paper we
get rid of the small factor $2^{O(\alpha^2(n))}$, and obtain the upper
bound $O(n^{3/2})$. Admittedly, the improvement is not large, and
achieves only a slight narrowing of the gap from the best known lower
bound, which is $\Omega(n^{4/3}\log\log n)$~\cite{E60}, but is
nevertheless the first improvement of the bound of \cite{CEGSW}, more
than 20 years after its establishment.

The proof of the new bound is based on the recently introduced
\emph{polynomial partitioning} technique of Guth and Katz~\cite{GK2}
(also see Kaplan et al.~\cite{KMS} for an expository introduction).
An additional goal of the present
paper is to highlight certain technical issues (specifically,
multi-level polynomial partitions) that might arise in the application
of the new approach. These issues are relatively simple to handle
for the problem at hand, but treating them in full generality
is still an open issue.

\heading{Zahl's work. }
After we finished a draft of this paper, in early 2011,
we learned that Zahl \cite{Za}
had independently obtained the same bound on unit distances in $\R^3$
(and, actually, a more general result concerning incidences
of points with suitable surfaces in $\R^3$), using the same general
approach. We believe that our treatment may still be of some
interest: First, the details of our arguments differ from those
of Zahl at some points, and since the general problem of the multi-level
decomposition alluded to above remains
unresolved, even slight differences in the approaches may
become important in attacking the general question. Second,
our treatment is more pedestrian and assumes less background
in algebraic geometry than Zahl's, and thus it may be more
accessible for the community at large of
researchers in discrete geometry. So, while we respect
(and clearly acknowledge) the priority of Zahl's preprint,
we have still decided to make our manuscript publicly available.

\section{Analysis}

Let $P$ be a set of $n$ points in $\R^3$. For each $p\in P$ let
$\sigma_p$ denote the unit sphere centered at $p$, and let $\Sigma$
denote the collection of these spheres. Clearly, the number
of unit distances between pairs of points of $P$ is half the number of
incidences $I(P,\Sigma)$ of the points of $P$ with the spheres of
$\Sigma$. Our main result is the following theorem.
%%%%%%%%%%%%%%%%%%%%%%%%%%%%
\begin{theorem} \label{th:main}
$I(P,\Sigma) = O(n^{3/2})$.
In particular, the number of unit distances in any set of $n$ points
in $\reals^3$ is $O(n^{3/2})$.
\end{theorem}
%%%%%%%%%%%%%%%%%%%%%%%%%%%%

We first review the main algebraic ingredient of the analysis.

\heading{Polynomial partitions: A quick review.}
For the sake of completeness, and also for the second partitioning
step in our analysis, we provide a brief review of the
polynomial partitioning
technique of Guth and Katz~\cite{GK2}; see also \cite{KMS}. This
technique is based on the \emph{polynomial ham sandwich} theorem
of Stone and Tukey~\cite{ST}. Specifically, fix an integer $\Deg$
and put $M={\Deg+3\choose 3}-1$. Let $U_1,\ldots,U_M$ be $M$
arbitrary finite point sets in $\reals^3$ (the theorem holds
for more general sets, and in any dimension, but this version
suffices for our purposes). Apply the \emph{Veronese map}
$\varphi:\reals^3\mapsto\reals^M$, which maps a point
$(x,y,z)\in\reals^3$ to the $M$-tuple of the values at $(x,y,z)$
of all the $M$ nonconstant trivariate
monomials of degree at most $\Deg$.
Consider the images $\varphi(U_1),\ldots,\varphi(U_M)$ of our sets,
and apply the standard ham sandwich theorem
(see \cite{ST} and \cite[Chapter 3]{Ma}) to
these $M$ sets in $\reals^M$, to obtain a hyperplane $h$ that
\emph{bisects} each set $U_i$, in the sense that, for each
$i=1,\ldots,M$, at most $|U_i|/2$ points of $U_i$ lie on one
side of $h$ and at most $|U_i|/2$ points lie on the other side
(the remaining points of $U_i$ lie on $h$; their number can be
anything between $0$ and $|U_i|$). Consider the trivariate polynomial
$f = h\circ\varphi$ (here $h=0$ is the linear equation of our
hyperplane). Then $f$ is a trivariate polynomial (a linear 
combination of monomials) of degree at most $\Deg$ that bisects 
each of the sets $U_1,\ldots,U_M$, in the sense that, for each $i$,
$$
|U_i\cap\{f>0\}|,\; |U_i\cap\{f<0\}| \le |U_i|/2 .
$$
Guth and Katz construct a sequence of such polynomial ham sandwich
cuts, to partition the given point set $P$ into a
specified number $t$ of subsets, each consisting of at most $O(n/t)$
points. This technique first bisects the original point set $P$
into two halves, using a polynomial $f_1$. It then bisects each
of these two sets into two halves, using a second polynomial $f_2$,
bisects each of the four resulting subsets using a third polynomial
$f_3$, and so on, until the desired number $t$ of subsets is obtained
(or exceeded). The product $f=f_1f_2f_3\cdots$ of these bisecting
polynomials is the desired \emph{partitioning polynomial}, and, as
shown in \cite{GK2,KMS}, its degree is $\Deg=O(t^{1/3})$.
The resulting partition is not exhaustive, as some points of $P$ may
lie in the zero set $Z(f)$ of $f$.
Note that in general it makes sense to take $t\le n$.
If $t>n$ we can, following the technique used in \cite{EKS,GK},
find a polynomial $f$ of degree $O(n^{1/3}) = O(t^{1/3})$
that vanishes at all the points of $P$. In this case all the subsets
in the resulting partition of $P$ are empty, except for
$P\cap Z(f) = P$.

\medskip

\noindent{\bf First partition.}
For the proof of Theorem~\ref{th:main}, we set $t=n^{3/4}$,
so the degree of the resulting partitioning
polynomial $f$ is $\Deg=O(n^{1/4})$.
Denote the resulting subsets of the above partition of $P$ by
$P_1,\ldots,P_t$. Each of these subsets is of size at most $O(n/t)$,
and we also have a remainder subset $P_0$, contained in the zero set
$Z=Z(f)$ of $f$. We may also assume that, for each $i\ge 1$, $P_i$
lies in a distinct connected component $C_i$ of
$\reals^3\setminus Z$. (More specifically, the construction of Guth
and Katz ensures that no connected component of $\reals^3\setminus Z$
meets more than one of the sets $P_i$. Some of these sets, though, may
lie in several components, in which case we further split each
such set into subsets, one for each component that it meets.)

We note that the degree $\Deg$ could conceivably be much smaller. For
example, if $P$, or most of it, lies on an algebraic surface of small
degree (say, a plane or a quadric) then $f$ could be the polynomial
defining that surface, resulting in a trivial partitioning in which
all or most of the points of $P$ belong to $P_0$ and the degree of $f$
is very small. This potential variability of $\Deg$ will enter the
analysis later on.

We first bound the number of incidences between $P\setminus P_0$
and $\Sigma$. For this, we need to show that no sphere
crosses too many cells of the partition (that is, components
of $\R^3\setminus Z(f)$).  This can be argued as follows.

Let us fix a sphere $\sigma=\sigma_a\in\Sigma$.
The number of cells $C_i$ crossed by $\sigma$
is bounded from above by the number of components
of $\sigma\setminus Z(f)$.  

For bounding the latter quantity, as well as in some
arguments in the sequel, it is technically convenient
to use a rational parameterization of $\sigma$.
Specifically, we let
$\psi\colon\R^2\to\R^3$ be the \emph{inverse stereographic
projection} given by $\psi(u,v)=(\psi_x(u,v),\psi_y(u,v),
\psi_z(u,v))$, where
$$
\psi_x(u,v) = x_0 + \frac{2u}{u^2+v^2+1},\quad\quad
\psi_y(u,v) = y_0 + \frac{2v}{u^2+v^2+1},\quad\quad
\psi_z(u,v) = z_0 + \frac{u^2+v^2-1}{u^2+v^2+1},
$$
and $(x_0,y_0,z_0)$ is the center of $\sigma$. Then
$\psi$ is a homeomorphism between the $uv$-plane and
the sphere $\sigma$ ``punctured'' at its north pole. This missing
point will not affect our analysis if we choose a generic coordinate
frame, in which no pair of points of $P$ are co-vertical. (Since the
center of each ball is a point in $P$, no point can reside at the
north pole of a ball in such a generic coordinate frame.)

Let us consider the composition $f\circ\psi$
(i.e., $f\circ\psi(u,v)=f(\psi_x(u,v),\psi_y(u,v),\psi_z(u,v))$);
this is a rational function, which we can write as a
quotient $\frac {f^*(u,v)}{q(u,v)}$ of two polynomials
(with no common factor). For analyzing the zero set,
it suffices to consider the
numerator $f^*(u,v)$, which is a polynomial of degree  $O(\Deg)$.

If $f^*$ vanishes identically then $\sigma\subset Z(f)$ and thus
$\sigma$ does not cross any cell $C_i$ of the partition.
Otherwise, the number of components of $\sigma\setminus Z(f)$
is no larger than the number of components of $\R^2\setminus Z(f^*)$,
and for these, we use the $d=2$ case of the following result.
%variant of 
%Harnack's theorem:\footnote{Harnack's theorem itself deals
%with the number of components of $Z(f)$, while we
%need to consider the components of the complement.}

\begin{lemma}\label{l:cpts}
Let $f$ be a real polynomial of degree $D$ in $d$ variables. Then
the number of connected components of $\R^d\setminus Z(f)$ is at
most $6(2D)^d$.
\end{lemma}

This lemma follows, for example, from Warren \cite[Theorem~2]{Wa}
(also see \cite{BPR} for an exposition, and \cite{Akama-al} for a
neatly simplified proof).

From Lemma \ref{l:cpts} we get that the number of connected
components of $\R^2\setminus Z(f^*)$ is $O(\dg(f^*)^2) = O(\Deg^2)$.
We thus conclude that each sphere $\sigma=\sigma_a\in\Sigma$ crosses
at most $O(\Deg^2) = O(n^{1/2})$ cells $C_i$ of the partition.

Hence the overall number of sphere-cell crossings is
$O(n\Deg^2) = O(n^{3/2})$. Situations in which a sphere
$\sigma$ crosses a cell $C_i$ and is incident to at most
two points of the corresponding subset $P_i = P\cap C_i$
generate a total of $O(n^{3/2})$ incidences within this subset.
Otherwise, for a cell $C_i$, its corresponding subset
$P_i$, and a fixed point $p\in P_i$, the number of
spheres that are incident to $p$ and contain at least two other
points of $P_i$, is at most $2{|P_i|-1\choose 2} \le |P_i|^2$,
because any pair of points $q,r\in P_i\setminus\{p\}$ determine
at most two unit spheres that are incident to $p,q,r$. Hence the
number of incidences of the points of $P_i$ with spheres that are
incident to at least three points of $P_i$ is at most
$|P_i|^3\le (n/t)^3 = O(n^{3/4})$. Summing over all subsets
$P_i$, we get a total of $O(n^{3/2})$ such incidences.\footnote{%
  Alternatively, we can use the K\H{o}v\'ari--S\'os--Tur\'an theorem
  (see \cite{PA}) on the maximum number of edges in a bipartite graph
  with a forbidden $K_{r,s}$ subgraph, as was done in many previous
  papers; this comment applies to several similar arguments 
  in the sequel.}

\noindent{\bf Remark.} 
(Although the full significance of this remark will become clearer
later on, we nevertheless make it early in the game.)
There are well-known papers
in real algebraic geometry estimating the number of components
of algebraic varieties in $\R^d$, or more generally,
the complexity of an arrangement of zero sets of
polynomials in $\R^d$ (Oleinik and Petrovski\v{\i}, Milnor, Thom,
and Warren---see, e.g., \cite{BPR} for references). 
In the arguments used so far, and also in the sequel, 
we need bounds in a somewhat different
setting, namely, when the arrangement is not in $\R^d$, but
within some algebraic variety. This setting was considered
by Basu, Pollack, and Roy \cite{bpr-ncdfp}; however, their bound is not
sufficiently sharp for us either, since it assumes the same upper bound
both on the degree of the polynomials defining the arrangement
and those defining the variety. Prompted by our question,
Barone and Basu~\cite{BB} proved a bound in this setting involving
two degree parameters: they consider a $k$-dimensional variety $V$
in $\R^d$ defined by polynomials of degree at most $D$, and an
arrangement of $n$ zero sets of polynomials of degree at most $E$ within $V$,
and they bound the number of cells, of all dimensions,
in the arrangement by $O(1)^d D^{d-k} (nE)^k$.
A weaker bound of a similar kind was also derived independently
by Solymosi and Tao \cite[Theorem~B.2]{SoTa}.
However, in our proof, we will eventually need three different
degree parameters (involving spheres intersecting a variety defined by
two polynomials of two potentially different degrees; in this case one
of the degrees is $2$, the degree of the polynomial equation of a
sphere), and thus we cannot refer to \cite{BB,SoTa} directly. 
We provide elementary ad-hoc arguments
instead (aimed mainly at readers not familiar with the
techniques employed in \cite{BB,SoTa}). If the multi-level
polynomial partition method should be used in dimensions higher
than 3, a more systematic approach will be needed to bound
the appropriate number of components. We believe that the
approach of  \cite{BB} should generalize to an arbitrary
number of different degree parameters, but there are
several other obstacles to be overcome along the way; see
Section~\ref{sec:conc} for a discusssion.

\noindent{\bf\boldmath Bounding $I(P_0,\Sigma)$.} It therefore remains to
bound $I(P_0,\Sigma)$. Here is an informal overview of this second
step of the analysis. We apply the polynomial partitioning procedure
to $P_0$, using a second polynomial $g$ (which again is the product
of logarithmically many bisecting polynomials). For a good choice of
$g$, we will obtain various subsets of $P_0$ of roughly equal sizes,
lying in distinct components of $Z(f)\setminus Z(g)$, and a
remainder subset $P_{00}\subset Z(f)\cap Z(g)$. Again, for a good
choice of $g$, $Z(f)\cap Z(g)$ will be a 1-dimensional curve, and it
will be reasonably easy to bound $I(P_{00},\Sigma)$. The situation
that we want to avoid is one in which $f$ and $g$ have a common
factor, whose 2-dimensional zero set contains most of $P_0$, in
which case the dimension reduction that we are after (from a
2-dimensional surface to a 1-dimensional curve) will not work.

To overcome this potential problem, we first factor $f$ into
irreducible factors $f=f_1f_2\cdots f_r$ (recall that in the
construction of \cite{GK2} $f$ is the product of logarithmically
many factors, some of which may themselves be reducible). Denote
the degree of $f_i$ by $\Deg_i$, so $\sum_i \Deg_i = \Deg$.
By removing repeated factors from $f$, if any exist, we may assume
that $f$ is square-free; this does not affect the partition induced
by $f$, nor its zero set. Put
\begin{align*}
P_{01} &= P_{0}\cap Z(f_1) \\
P_{02} &= \left(P_{0}\setminus P_{01}\right) \cap Z(f_2) \\
&\cdots \\
P_{0i} &= \biggl(P_{0}\setminus \bigcup_{j<i} P_{0j}\biggr) \cap Z(f_i) \\
&\cdots
\end{align*}
This is a partition of $P_{0}$ into $r$ pairwise disjoint subsets.
Put $m_i = |P_{0i}|$ for $i=1,\ldots,r$; thus, $\sum_i m_i \le n$. We
will bound $I(P_{0i},\Sigma)$ for each $i$ separately and then add up
the resulting bounds to get the desired bound on $I(P_0,\Sigma)$.

\noindent{\bf Second partition. } We will bound the
number of incidences between $P_{0i}$ and $\Sigma$ using the following
lemma, which is the core of (this step of) our analysis. 

\begin{lemma} \label{lem:2nd}
Let $f$ be an irreducible trivariate polynomial of degree $\Deg$, let
$Q$ be a set of $m$ points contained in $Z(f)$, and let $\Sigma$
be a set of $n\ge m$ unit spheres in $\reals^3$. Then
$$
I(Q,\Sigma) = O\Bigl(
m^{3/5}n^{4/5}\Deg^{2/5} + n\Deg^2 \Bigr) .
$$
\end{lemma}

\noindent{\bf Remark.} 
When $D=1$ (all the points of $Q$ are co-planar), the bound in the
lemma becomes $O(m^{3/5}n^{4/5} + n)$, a special case (when $m\le n$)
of the bound $O(m^{3/5}n^{4/5} + n + m)$, which is a well known upper
bound on the number of incidences between $m$ points and $n$ circles
in the plane (see, e.g., \cite{CEGSW,PS}). 
%\micha{A better ref, without the old $\delta$'s in the exponents?}
In our case, the circles are
the intersections of the spheres of $\Sigma$ with the plane (where
each resulting circle has multiplicity at most $2$).

The main technical step in proving Lemma~\ref{lem:2nd} is
expressed by the following lemma.

\begin{lemma} \label{lem:fg}
Given an irreducible trivariate polynomial $f$ of degree $\Deg$, a
parameter $\Degg\ge\Deg$, and a finite point set $Q$ in $\R^3$, there is a
polynomial $g$ of degree at most $\Degg$, co-prime with $f$, which
partitions $Q$ into subsets $Q_0\subseteq Z(g)$ and
$Q_1,\ldots,Q_t$, for $t=\Theta(\Deg\Degg^2)$, so that each $Q_i$,
for $i=1,\ldots,t$, lies in a distinct component of $\R^3\setminus
Z(g)$, and $|Q_i| =O(|Q|/t)$.
\end{lemma}

Note the similarity of this lemma to the standard polynomial
partitioning result, as used in the first partitioning step. The
difference is that, to ensure that $g$ be co-prime with $f$, we
pay the price of having only $\Theta(\Deg\Degg^2)$ parts in the
resulting partition, instead of $\Theta(\Degg^3)$.

\noindent{\bf Proof of Lemma \ref{lem:fg}.} As in the standard
polynomial partitioning technique, we obtain $g$ as the product of
logarithmically many bisecting polynomials, each obtained by
applying a variant of the polynomial ham sandwich theorem to a
current collection of subsets of $Q$. The difference, though, is
that we want to ensure that each of the bisecting polynomials is not
divisible by $f$; since $f$ is irreducible, this ensures
co-primality of $g$ with $f$. Reviewing the construction of
polynomial ham sandwich cuts, as outlined above, we see that all
that is needed is to come up with some sufficiently large finite 
set of monomials, of an appropriate maximum degree,
so that no nontrivial linear combination of
these monomials can be divisible by $f$. We then use a restriction
of the Veronese map defined by this subset of monomials, and the
standard ham-sandwich theorem in the resulting high-dimensional space,
to obtain the desired polynomial.

Let $x^iy^jz^k$ be the \emph{leading term} of $f$, in the sense that
$i+j+k=\Deg$ and $(i,j,k)$ is largest in the lexicographical order
among all the triples of exponents of the monomials of $f$  of
degree $\Deg$.
Let $s$ be the desired number of sets that we want a single
partitioning polynomial to bisect. For that we need a space of $s$
monomials whose degrees are not too large and which 
span only polynomials not divisible by $f$. If, say,
$s < \left(\frac{\Deg}{3} \right)^3$ then we can use all monomials
$x^iy^jz^k$ such that $i,j,k \le s^{1/3} < \Deg/3$. Clearly, any 
nontrivial linear
combination of these monomials cannot be divisible by $f$. In this
case the degree of the resulting partitioning polynomial is
$\Theta(s^{1/3})$.  If $s > \left(\frac{\Deg}{3} \right)^3$ then 
we take the set of all monomials
$x^{i'}y^{j'}z^{k'}$ that satisfy $i'<i$ or $j'<j$ or $k'<k$, and
$\max\{i',j',k'\} \le \hat{\Deg}$ for  a suitable integer $\hat{\Deg}$, 
which we specify below (the actual degree of the bisecting polynomial 
under construction will then be at most $3\hat{\Deg}$). 
Any nontrivial polynomial $h$ which is a linear combination of
these monomials cannot be divisible by $f$. Indeed, if $h=fh_1$ for
some polynomial $h_1$ then the product of the leading terms of $f$
and of $h_1$ cannot be canceled out by the other monomials of the
product, and, by construction, $h$ cannot contain this monomial. The
number of monomials in this set is
$\Theta(i\hat{\Deg}^2+j\hat{\Deg}^2+k\hat{\Deg}^2) =
\Theta(\Deg\hat{\Deg}^2)$. We thus pick $\hat{\Deg}=
\Theta(({s}/{\Deg})^{1/2})$ so that we indeed get $s$ monomials. 
As noted above, the degree of the resulting bisecting polynomial 
in this case is $O((s/{\Deg})^{1/2})$.

We now proceed to construct the required partitioning of $Q$ into $t$ sets,
by a sequence of about $\log t$ polynomials $g_0,g_1,\ldots$, where $g_j$
bisects $2^{j}$ subsets of $Q$, each of size at most $|Q|/2^{j}$. For
every $j$ such that $s = 2^{j} < \left(\frac{\Deg}{3}\right)^3$ we
construct, as shown above, a polynomial of degree  $O(s^{1/3}) =
O(2^{j/3})$. For the indices $j$ with $s = 2^{j} >
\left(\frac{\Deg}{3}\right)^3$ we construct a polynomial of degree
$O(\left({s}/{\Deg}\right)^{1/2}) = O({2^{j/2}}/{\Deg^{1/2}})$. Since the
upper bounds on the degrees of the partitioning polynomials increase
exponentially with $j$, and since the number of parts that we want is
$\Omega(\Deg^3)$, it follows that the degree of the product of the
sequence is $O(\left(t/\Deg\right)^{1/2})$. If we require this degree 
bound to be no larger than $\Degg$ then it follows that the size of the
partition that we get is $t=\Theta(\Deg\Degg^2)$. Clearly, $f$ does not
divide the product $g$ of the polynomials $g_j$, so $g$ satisfies all
the properties asserted in the lemma. 
$\Box$

\medskip

\noindent{\bf Remarks.}
(1) The analysis given above can be interpreted as being applied to
the \emph{quotient ring} $Q=\reals[x,y,z]/I$, where
$I=\langle f\rangle$ is the ideal generated by $f$.
General quotient rings are described in detail in, e.g.,
\cite{CLO1,CLO2}, but the special case where $I$ is generated by
a single polynomial is much simpler, and can be handled in the simple
manner described above, bypassing (or rather simplifying considerably)
the general machinery of quotient rings.
As a matter of fact, an appropriate extension of Lemma~\ref{lem:fg} to 
quotient rings defined by two or more polynomials is still an open
issue; see Section~\ref{sec:conc}.

\noindent 
(2) 
The set $Q$ is in fact contained in $Z(f)$, and the subset $Q_0$ 
is contained in $Z(f)\cap Z(g)$. 
However, except for the effect of this property
on the specific choice of
monomials for $g$, the construction considers $Q$ as an arbitrary
set of points in $\R^3$, and does not exploit the fact that
$Q\subset Z(f)$.

\medskip

\noindent{\bf Back to the proof of Lemma~\ref{lem:2nd}. }
We apply Lemma~\ref{lem:fg} to $Q$, now assumed to be contained in
$Z(f)$, and obtain the desired partitioning polynomial $g$.
We now proceed, based on the resulting partition of $Q$, to bound
$I(Q,\Sigma)$.

We need the following technical lemma, a variant of which has been
established and exploited in \cite{GK} and in \cite{EKS}. For the
sake of completeness we include a brief sketch of its proof, and
refer the reader to the aforementioned papers for further details.

\begin{lemma} \label{lem:x} 
{\rm (a)} Let $f$ and $g$ be two trivariate polynomials of respective
degrees $\Deg$ and $\Degg$. Let $\Pi$ be an infinite collection of parallel
planes such that, for each $\pi\in \Pi$, the restrictions of $f$ and
$g$ to $\Pi$ have more than $\Deg\Degg$ common roots. Then $f$ and $g$
have a (nonconstant) common factor. 

{\rm (b)} Let $f$ and $g$ be as in {\rm (a)}. If the intersection 
$Z(f)\cap Z(g)$ of their zero
sets contains a $2$-dimensional surface patch then $f$
and $g$ have a (nonconstant) common factor.
\end{lemma}

\noindent{\bf Proof sketch.} 
(a) Assume without loss of generality
that the planes in $\Pi$ are horizontal and that, if the number of
common roots in a plane is finite then these roots have 
different $x$-coordinates; both assumptions can be enforced by
an appropriate rotation of the coordinate frame. 
Consider the $y$-resultant
$r(x,z)={\rm Res}_y(f(x,y,z),g(x,y,z))$ of $f(x,y,z)$ and $g(x,y,z)$. This
is a polynomial in $x$ and $z$ of degree at most $\Deg\Degg$. If the plane
$z=c$  contains more than $\Deg\Degg$ common roots then $r(x,c)$, which is
a polynomial in $x$, has more than $\Deg\Degg$ roots, and therefore 
it must be identically zero. 
It follows that $r(x,z)$ is identically zero on
infinitely many planes $z=c$, and therefore, it must be identically zero.
(Its restriction to an arbitrary non-horizontal line $\ell$ has infinitely 
many roots and therefore it must be identically zero on  $\ell$.) It
follows that $f(x,y,z)$ and $g(x,y,z)$ have a common factor (see
\cite[Proposition~1, page~163]{CLO1}).

(b) This follows from (a), since if $Z(f)$ and $Z(g)$ contain a
2-dimensional surface patch, then they must have infinitely many zeros on
infinitely many parallel planes. $\Box$

\noindent{\bf\boldmath Incidences outside $Z(g)$.} 
To prove Lemma~\ref{lem:2nd}, we first bound the number
of incidences of the points of a fixed subset $Q_j$, for $j\ge 1$,
with $\Sigma$, using the same approach as in the first partition.
That is, let $n_j$ denote the number of spheres of $\Sigma$ that
cross the corresponding cell $C_j$ {\em effectively}, in the sense
that $\sigma\cap Q_j\ne\emptyset$. Then we have $O(n_j)$ incidences
of the points of $Q_j$ with spheres that are incident to at most two
points of $Q_j$, and $O((m/t)^3)$ incidences with spheres that are
incident to at least three points. Summing over all sets, we get
\begin{equation} \label{eq:nj}
\sum_{j=1}^t I(Q_j,\Sigma) = O\biggl(m^3/t^2 + \sum_{j=1}^t n_j
\biggr).
\end{equation}
We estimate $\sum_j n_j$ by bounding the number of cells $C_j$
that a single sphere $\sigma\in\Sigma$ can cross effectively,
which we do as follows.

Take the same rational parametrization $\psi$
of $\sigma$ used in the analysis of the first
partitioning step.
 Let $f^*(u,v)$ and $g^*(u,v)$ be the
polynomials obtained from $f\circ \psi$ and
$g\circ\psi$ by removing the common denominator of
these  rational functions. The  degrees of $f^*$ and $g^*$ are
 $O(\Deg)$ and $O(\Degg)$, respectively.

If $f^*$ vanishes identically on the $uv$-plane, then
$\sigma\subseteq Z(f)$; this is an easy situation that we will handle later
on. Otherwise, $Z(f^*)=\psi^{-1}(\sigma\cap
Z(f))$ is a 1-dimensional curve $\gamma$ in the $uv$-plane
(possibly degenerate, e.g., empty or consisting of isolated points), 
and $Q\cap\sigma$ is contained in $\psi(\gamma)$.

By construction, the number of cells $C_j$ that $\sigma$ crosses
effectively (so that it is incident to points of $Q_j$) is no larger
than the number of components of $Z(f^*)\setminus Z(g^*)$. This is
because each such cell $C_j$ contains at least one connected
component of $\psi(Z(f^*) \setminus Z(g^*))$. 

Now each component of $Z(f^*)\setminus Z(g^*)$
is either a full component of $Z(f^*)$, or a relatively
open connected portion of $Z(f^*)$ whose closure meets $Z(g^*)$.

Since $f^*$ is a bivariate polynomial, Harnack's
theorem~\cite{Harnack} asserts that the number of (arcwise)
connected components of $Z(f^*)$ is at most $1+{\dg(f^*)-1\choose 2}
= O(\Deg^2)$. 

For the other kind of components, choose a generic
sufficiently small value $\eps>0$, so that $f^*$ and $g^*\pm\eps$ do
not have a common factor.\footnote{Indeed, assuming that
$f^*$ and $g^*+\eps$ had a nonconstant common factor for infinitely
many values of $\eps$, then the same factor would occur for
two distinct values $\eps_1$ and $\eps_2$ of $\eps$,
and thus it would have to divide $\eps_1-\eps_2$, which is impossible.}
Then each component of $Z(f^*)\setminus
Z(g^*)$ of the second kind must contain a point at which 
$g^*+\eps=0$ or $g^*-\eps=0$. Hence, the number of such components
is at most the number of such common roots, which, by B\'ezout's
theorem (see, e.g.,~\cite{CLO2}) 
is\footnote{The $O(\Deg\Degg)$ bound for the
number of components of $Z(f^*)\setminus Z(g^*)$ is also a direct
consequence of the main result of Barone and Basu~\cite{BB}.}
$O(\dg(f^*)\dg(g^*)) = O(\Deg\Degg)$.

Since $\Degg\ge\Deg$, we conclude that the number of
cells $C_j$ crossed effectively by $\sigma$ is $O(\Deg\Degg)$, which
in turn implies that $\sum_j n_j = O(n\Deg\Degg)$. Substituting this
in (\ref{eq:nj}) and recalling that $t=\Theta(\Deg\Degg^2)$, we get
\begin{equation} \label{eq:inc1}
\sum_{j=1}^r I(Q_j,\Sigma) = O\biggl(\frac{m^3}{\Deg^2\Degg^4} +
n\Deg\Degg \biggr).
\end{equation}

We have left aside the case where $\sigma\subseteq Z(f)$.
Since $f$ is irreducible, and so is $\sigma$, we must have
$\sigma=Z(f)$ in this case (recall Lemma~\ref{lem:x}(b)).
The analysis proceeds as
above for every sphere $\sigma' \not= \sigma$, and the number of
incidences with $\sigma$ itself is at most $m$, a bound that is
subsumed by the bound asserted in the lemma (recall that $m\le n$).

We note that in the ongoing analysis $D$ is the actual degree of the
irreducible factor of $f$ under consideration, but $E$ is only  a
chosen upper bound for $\dg(g)$, whose actual value may be smaller
(as may have been the case with $f$).

To optimize the bound in (\ref{eq:inc1}), we choose
\begin{equation} \label{eq:E}
\Degg = \max\Bigl\{ \frac{m^{3/5}}{n^{1/5}\Deg^{3/5}},\;\Deg\Bigr\} ,
\end{equation}
and observe that the first term dominates when $\Deg \le
m^{3/8}/n^{1/8}$. Assuming that this is indeed the case, we get
\begin{equation} \label{eqn:1x}
\sum_{j} I(Q_j,\Sigma) = O(m^{3/5}n^{4/5}D^{2/5}) \ .
\end{equation}
If $\Deg > m^{3/8}/n^{1/8}$ then we have $\Degg = \Deg$, and the bound
(\ref{eq:inc1}) becomes
\begin{equation} \label{eq:inc2}
\sum_j I(Q_j,\Sigma) =
O\left(\frac{m^3}{\Deg^6} + n\Deg^2 \right) = O\left(n\Deg^2 \right) .
\end{equation}
Thus, $I(Q\setminus Q_0,\Sigma)$ satisfies the bound asserted in
the lemma, and it remains to bound $I(Q_0,\Sigma)$.

\medskip

\noindent{\bf \boldmath Incidences within $Z(f)\cap Z(g)$.}
Recall that $Q_0$ is contained in the curve
$\delta=Z(f)\cap Z(g)$, which by Lemma \ref{lem:x}(b) is (at most)
1-dimensional. 

Fix a sphere $\sigma\in\Sigma$ that does not
coincide with $Z(f)$, let $\psi_\sigma$ be the corresponding
rational parameterization of $\sigma$,
and let $f_\sigma^*$ and $g_\sigma^*$ be as defined in the
preceding analysis. 

If $g_\sigma^*$ is identically 0, then we have
$\sigma\subseteq Z(g)$, and the irreducible polynomial
defining $\sigma$ is a factor of $g$ by Lemma~\ref{lem:x}. 
Thus, the number of such $\sigma$'s is $O(\Degg)$, and together they
can contribute at most $O(m\Degg)$ incidences,
which is bounded from above by the right-hand side of~(\ref{eq:inc1}).

Now we assume that both $f_\sigma^*$ and $g_\sigma^*$ are nonzero,  we let
let $h_\sigma^*$ denote the greatest common divisor of
$f_\sigma^*$ and $g_\sigma^*$, and put
$f_\sigma^*=f_{1\sigma}^*h_\sigma^*$ and
$g_\sigma^*=g_{1\sigma}^*h_\sigma^*$. Then
$\psi_\sigma^{-1}(\sigma\cap\delta)$ is the union of $Z(h_\sigma^*)$ and of
$Z(f_{1\sigma}^*)\cap Z(g_{1\sigma}^*)$.  
Using B\'ezout's theorem as above, we have
$|Z(f_{1\sigma}^*)\cap Z(g_{1\sigma}^*)| = O(\Deg\Degg)$;  
summing this bound over
all spheres, we get at most $O(n\Deg\Degg)$ incidences, a bound
already subsumed by (\ref{eq:inc1}).

It remains to account for incidences of the following
kind (call them \emph{$h^*$-incidences}): 
a point $q\in Q_0\cap\sigma$ lying in $\psi_\sigma(Z(h_\sigma^*))$.
Let us call such a point $q$ \emph{isolated in $\sigma$}
if it is an isolated point of $\psi_\sigma(Z(h_\sigma^*))$;
i.e., there is a neighborhood of $q$ in $\sigma$ intersecting
$\psi_\sigma(Z(h_\sigma^*))$ only at $q$.

The homeomorphism $\psi_\sigma^{-1}$ maps the isolated points $q$
on $\sigma$ to isolated
zeros of $h_\sigma^*$ in the $uv$-plane, in a one-to-one fashion.
Since $\deg(h^*_\sigma)=O(\Deg)$, $Z(h_\sigma^*)$
has at most $O(\Deg^2)$ components (Harnack's theorem again),
and thus the overall number of isolated incidences is $O(n\Deg^2)$.

Finally, to account for non-isolated $h^*$-incidences,
let us fix a point $q\in Q_0$,
and consider the collection $\hat{\Sigma}_q$ consisting of all
spheres $\sigma\in\Sigma$ that contain $q$ such that
$q$ forms a non-isolated $h^*$-incidence with $\sigma$.
We claim that $|\hat{\Sigma}_q|=O(\Deg\Degg)$.

For $\sigma\in \hat{\Sigma}_q$, the set $\psi_\sigma(Z(h_\sigma^*))$
contains a curve segment $\beta_{q,\sigma}$ ending at $q$.
Let us call $\beta_{q,\sigma}$ and $\beta_{q,\sigma'}$
\emph{equivalent} if they coincide in some neighborhood of $q$.
If $\beta_{q,\sigma}$ and $\beta_{q,\sigma'}$ are not equivalent,
then in a sufficiently small neighborhood of $q$
they intersect only at $q$ (since they are arcs of algebraic curves). 

We also note that a given $\beta_{q,\sigma}$ can be equivalent to
$\beta_{q,\sigma'}$ for at most one $\sigma'\ne\sigma$; this
is because the common portion $\beta_{q,\sigma}\cap \beta_{q,\sigma'}$ 
of the considered curve segments has to be
contained in the intersection circle $\sigma\cap\sigma'$,
and that circle intersects any other sphere $\sigma''\in\Sigma$
in at most two points. Thus, $|\hat\Sigma_q|$ is at most twice
the number of equivalence classes of the curve segments~$\beta_{q,\sigma}$.

Let us fix an auxiliary sphere $S$ of a sufficiently small radius $\rho$
around $q$, so that each $\beta_{q,\sigma}$ intersects $S$
at some point $x_{\sigma}$. Let $S'$ be the sphere around $q$
of radius $\rho/10$, say; we choose a point $y\in S'$
uniformly at random, and let $\pi$ be the plane tangent
to $S'$ at $y$. Then, for each $\sigma\in\hat\Sigma_{q}$,
$\pi$ separates $x_\sigma$ from $q$ with probability
at least $\frac 13$, say, and thus, by continuity, it intersects
$\beta_{q,\sigma}$. Hence there is a specific  $y_0\in S'$
such that the corresponding tangent plane
$\pi_0$ intersects $\beta_{q,\sigma}$ for at least
a third of the spheres $\sigma\in\hat \Sigma_q$.

Moreover, we can assume that such a $\pi_0$ intersects each
$\beta_{q,\sigma}$ in such a way that
all planes $\pi$ parallel to $\pi_0$ and sufficiently close
to it intersect $\beta_{q,\sigma}$ as well.
Then an application of Lemma \ref{lem:x}(a)
allows us to assume that the restrictions of $f$ and $g$ to 
some $\pi$ as above
are bivariate polynomials, with at most $\Deg\Degg$ common
roots. Hence $\pi$ intersects at most
$O(\Deg\Degg)$ of the curves $\beta_{q,\sigma}$, and 
so $|\hat\Sigma_q|=O(\Deg\Degg)$.

%Indeed, a generic plane $\pi$ passing sufficiently close to $p$ 
%will cross all the curves $\psi(Z(h_\sigma^*))$ for 
%$\sigma\in \hat{\Sigma}_p$. 
%(For this, we need to assume that $p$ is not a cusp
%of any of these curves. If such degeneracies arise, we can shoose a
%generic orientation of $\pi$, and then ensure that it crosses at
%least half of the curves if we place it at an appropriate side of
%$p$.) 

Altogether, we can bound the number
of $h^*$-incidences by  $O(n\Deg^2+m\Deg\Degg)$, which 
does not exceed  the earlier estimate $O(n\Deg\Degg)$. Hence
choosing $E$ as in (\ref{eq:E}), the incidences within $\delta$ do
not affect either of the asymptotic bounds (\ref{eqn:1x}),
(\ref{eq:inc2}).

This completes the proof of Lemma~\ref{lem:2nd}.  $\Box$

\medskip

\noindent{\bf Finishing the proof of Theorem~\ref{th:main}.}
We recall that in the first partitioning step, the set
 $P_0=P\cap Z(f)$ has been partitioned into the subsets 
$P_{01},\ldots,P_{0r}$.
Each $P_{0i}$ consists of $m_i$ points, and it is contained 
in $Z(f_i)$, where $f_i$ is an irreducible factor of $f$,
with $\deg(f_i)=\Deg_i$.
By Lemma~\ref{lem:2nd} we have
$$
\sum_{i=1}^r I(P_{0i},\Sigma)  =
O\biggl(\sum_{i=1}^r m_i^{3/5}n^{4/5}\Deg_i^{2/5}
+\sum_{i=1}^r n\Deg_i^2  \biggr).
$$
For the first term on the right-hand side we 
use H\"older's inequality\footnote{%
  H\"older's inequality asserts that
  $\sum x_iy_i\le (\sum |x_i|^p)^{1/p}(\sum |y_i|^q)^{1/q}$ 
  for positive $p,q$ satisfying 
  ${\displaystyle \frac1p + \frac1q = 1}$. 
  Here we use it with $p=\frac 53$, $q=\frac 52$, 
  $x_i=m_i^{3/5}$, and $y_i=D_i^{2/5}$.}
and the inequalities $\sum_{i=1}^r D_i\le D=O(n^{1/4})$ and
$\sum_{i=1}^r m_i\le n$. Thus,
$$
n^{4/5}\sum_{i=1}^r m_i^{3/5}\Deg_i^{2/5}\le 
n^{4/5} \biggl(\sum_i m_i\biggr)^{3/5} \biggl(\sum_i \Deg_i\biggr)^{2/5}
\le O(n^{4/5}n^{3/5}\Deg^{2/5}) = O(n^{3/2}).
$$
For the remaining term we have
$$
\sum_{i=1}^r
n\Deg_i^2 \le n\Deg\cdot  \sum_{i=1}^r \Deg_i  \le n\Deg^2 =
O(n^{3/2}).
$$
We thus get a total of $O(n^{3/2})$ incidences, thereby completing the
proof of the theorem.
$\Box$

\section{Discussion}
\label{sec:conc}

The main technical ingredient in the analysis, on top of the
standard polynomial partitioning technique of Guth and Katz, is the
recursion on the dimension of the ambient manifold containing the
points of $P$. This required a more careful construction of the
second partitioning polynomial $g$ to make sure that it is co-prime
with the first polynomial $f$. It is reasonably easy to perform the
first such recursive step, as done here and also independently by
Zahl~\cite{Za}, but successive recursive steps become trickier. In
such cases we have several co-prime polynomials, and we need to
construct, in the quotient ring of their ideal, a polynomial ham
sandwich cut of some specified maximum degree with sufficiently many
monomials. Such higher recursive steps will be needed when we
analyze incidences between points and surfaces in higher dimensions.
At the moment there does not seem to be an efficient procedure for
this task. Another recent paper where similar issues arise is by
Solymosi and Tao~\cite{SoTa}.

We also note that Zahl's study extends Theorem~\ref{th:main}
to incidences between points and more general surfaces in three
dimensions. The analysis in our study can also be similarly extended
(at the price of making some of the arguments more complicated),
but, since our goal had been to improve the bound on unit distances,
we have focused on the case of unit spheres.

\end{document}